\documentclass{amsart}
\usepackage{amssymb,amsmath,amscd,xy,graphicx,textcomp}
\newtheorem{theorem}{Theorem}[section]

\newtheorem{corollary}[theorem]{Corollary}
\newtheorem{definition}[theorem]{Definition}

\newtheorem{example}[theorem]{Example}
\newtheorem{proposition}[theorem]{Proposition}

\xyoption{arrow}

\xyoption{matrix}

\def\C{\mathbb{C}}
\def\R{\mathbb{R}}
\def\Z{\mathbb{Z}}

\def\br{\mathbf{r}}

\def\tree{\mathcal{T}}

\newcommand{\ZZ}{\mathbb Z}

\begin{document}

\title[Gorenstein algebras of weighted trees]{Gorenstein semigroup algebras of weighted trees and ordered points on the projective line}

\author{Christopher Manon}
\thanks{The author was supported by the NSF FRG grant DMS-0554254.}
\date{\today}

\begin{abstract}
We determine exactly which semigroup algebras of weighted trees are Gorenstein. These algebras arise as toric degenerations of projective coordinate rings of moduli of weighted points on the projective line. As a corollary, we find exactly when these families of algebras are Gorenstein as well. 
\end{abstract}

\maketitle

\tableofcontents

\smallskip

\section{Introduction}

Let $\mathfrak{k}$ be a field, then a $\ZZ$-graded $\mathfrak{k}$-algebra $R$ is said to be $Gorenstein$ if the Matlis Dual
$$H_{m}^{dim(R)}(R)^* = \underline{Hom}_{\mathfrak{k}}(H_{m}^{dim(R)}(R), \mathfrak{k}),$$
is isomorphic to grade-shifted copy $R(-a)$ of $R$. Here $m$ is the 
maximal ideal generated by elements in $R$  of positive degree, and $\underline{Hom}_{\mathfrak{k}}(-,-)$ 
is the functor of graded $\mathfrak{k}$-morphisms.  The number $a$ is called
the $a$-invariant of the graded Gorenstein algebra $R$.  We refer the reader to
the book by Bruns and Herzog, \cite{BH} for this and all other relevant definitions.
We study this property for a specific family of normal semigroup
algebras $\C[S_{\tree}(\br)].$  Here $\tree$ is a trivalent tree 
with $n$ ordered leaves, 
and $\br$ is an $n$-tuple of positive
integers which sums to an even number. We let $E(\tree)$ denote the 
set of edges of $\tree$, $L(\tree)$ denote the set 
of edges incident on a leaf of $\tree$,
$I(\tree)$ denote the set $E(\tree) \setminus L(\tree)$, 
and $V(\tree)$ denote the set of vertices of $\tree$. 
For each trinode $\tau \in V(\tree)$ we number the three edges incident on $\tau$
in some way with $\{1, 2, 3\}$, and we denote the $i$-th such edge by $(\tau, i).$  
When we are dealing with a convex cone or polytope $P$ we denote
the lattice points in the interior by $int(P)$. We will work over $\C,$ but our
analysis works for an arbitrary field $\mathfrak{k}$. 

\begin{definition}
The semigroup $S_{\tree}(\br)$ is graded, and the
$k$-th graded component is the set of
weightings $\omega: E(\tree) \to \ZZ_+$ defined by the following conditions. 

\begin{enumerate}
\item For all $\tau \in V(\tree)$ the numbers $\omega(\tau, i)$
satisfy $|\omega(\tau, 1) - \omega(\tau, 2)| \leq \omega(\tau, 3) \leq |\omega(\tau, 1) + \omega(\tau, 2)|.$
\item  $\Sigma_{i = 1}^3 \omega(\tau, i)$ is even. 
\item  For all $v_m \in L(\tree)$, $\omega(v_m) = k\br_m.$
\end{enumerate}
\end{definition}
The expressions in item $1$ above are called
the triangle inequalities, item $2$ is referred to 
as the parity condition. From now on, when three numbers $A$, $B$ and $C$ satisfy both
the parity condition and the triangle inequalities,
we write $\Delta_2(A, B, C).$  Notice that for a single trinode $\tau,$
these conditions are symmetric in the $(\tau, i).$
The graded components of the semigroups $S_{\tree}(\br)$ are subsets of the semigroup $S_{\tree}$
of non-negative integer weightings of $\tree$ which satisfy conditions
$(1)$ and $(2).$   This semigroup can be described as the lattice points in the 
convex cone $P_{\tree} \subset \R^{E(\tree)}$ defined by positive real vectors satisfying $(1),$
with respect to the lattice $L_2(\tree)$ of integer vectors satisfying $(2).$  When $\tree$ only has three leaves
we denote this cone by $P_3.$ The members of the first graded piece of $S_{\tree}(\br)$ are the lattice points in a cross-section $P_{\tree}(\br)$
of $P_{\tree},$ defined by specializing the weights in $L(\tree)$ to the entries in the vector $\br.$
\\

Presentations of these semigroups and their
associated algebras were constructed by the author
in \cite{M}.  In \cite{HMSV}, Howard, Millson, Snowden, and Vakil
independently constructed presentations of $\C[S_{\tree}(\br)]$
in order to find presentations of a projective coordinate ring
of the moduli space
of $\br$-weighted points on $\mathbb{P}^1$, denoted $M_{\br}$.  
The embedding they studied comes from the homeomorphism of projective varieties, 
$M_{\br} \cong Gr_2(\C^n) //_{\br} T$, 
where the right hand side is the $\br$-weight variety of $Gr_2(\C^n).$
They constructed a toric degeneration of this algebra, $\C[M_{\br}]$,
to $\C[S_{\tree}(\br)]$ for
each tree $\tree$ by means of the Speyer-Sturmfels 
\cite{SpSt} toric deformations of $Gr_2(\C^n),$
which deform the projective coordinate ring given by the Pl\"ucker 
embedding to $\C[P_{\tree}],$ the algebra of lattice points for the cone 
$P_{\tree}$, see \cite{HMSV} and \cite{SpSt} for details. 
Also, Sturmfels and Xu \cite{StXu} have shown the Cox ring 
$R_{n-1, n-3}$ of the blow up of $\mathbb{P}^{n-3}$
at $n-1$ points to be isomorphic (as a multigraded algebra) 
to the Pl\"ucker algebra of $Gr_2(\C^n).$ The multigrading, given
by the Picard group of this blow-up, corresponds to the multigrading on $\C[P_{\tree}]$
given by the weights on elements of $L(\tree).$  This then implies
that the subrings of $R_{n-1, n-3}$ given by taking invariants
with respect to a character of the corresponding ``Picard torus'' 
are isomorphic to $\C[M_{\br}].$  
The relevance of all the above work to this project 
comes from the following theorem.  

\begin{theorem}\label{deform}
The graded algebra 
$\C[M_{\br}]$ is Gorenstein if and only if
$\C[S_{\tree}(\br)]$ is Gorenstein, for any 
tree $\tree$.  
\end{theorem}
 \noindent
This follows from a theorem of Stanley, appearing as
corollary 4.3.8 in \cite{BH}, which characterizes
the Gorenstein properties of domains
by details of their Hilbert function.  We get the above theorem
because both $\C[M_{\br}]$ and $\C[S_{\tree}(\br)]$ are
domains, and have the same
Hilbert function by a standard theorem of deformation
theory. We are therefore able to study
the Gorenstein property for projective coordinate rings
of weighted points on $\mathbb{P}^1$, weight varieties of 
$Gr_2(\C^n)$, and the invariant subrings of $R){n-1, n-3}$ by learning about 
the family $\C[S_{\tree}(\br)]$ of normal semigroup algebras.
We also get the following corollary. 

\begin{corollary}\label{equiv}
$\C[S_{\tree}(\br)]$ is Gorenstein if and only
if $\C[S_{\tree'}(\br)]$ is Gorenstein for any trivalent trees
$\tree$, $\tree'$ with the same number of leaves.  
\end{corollary}
Since we are dealing with algebras generated
by the lattice points of convex cones we are
allowed the full range of available literature
on the commutative algebra of affine semigroup algebras, in particular
the following proposition, which is a consequence 
of corollary 6.3.8 in \cite{BH}.

\begin{proposition}
Let $S_P$ be the semigroup given by the lattice
points in a convex cone $P$.  Then the algebra
$\C[S_P]$ is Gorenstein if and only if 
there is a lattice point $\omega \in int(P)$ with $int(P) = \omega + P.$
Furthermore, in the presence of a grading, 
we have $a(\C[S_P]) = deg(\omega).$ 
\end{proposition}
\noindent
This proposition follows from the fact that 
the ideal $(int(P))\C[S_P]$ can identified
with the $canonical$ $module$ of the algebra $\C[S_P]$
(resp. the *-canonical module in 
the presence of a grading), see \cite{BH} for details.
We wish to prove this property for $S_{\tree}(\br)$, seen as the lattice points in the cone
over $P_{\tree}(\br)\times\{1\}$ in $\R^{I(\tree)} \times \R$
with respect to the product lattice $L_2 \times \Z$.  In order to do
this, we break the problem into two parts. First, we analyze when
some Minkowski sum of $P_{\tree}(\br)$ contains a unique interior lattice point, a necessary
but not sufficient condition for the Gorenstein property. 

\begin{theorem}\label{T1}
$P_{\tree}(\br)$ has a unique interior
point if and only if $\br = \vec{2} + \vec{R}$ where
$\vec{R}$ is of one of the following types. 
\begin{enumerate}
\item $R_i = \Sigma_{j \neq i} R_j$ for some $i$.
\item $\Delta_2(R_i, R_j, R_k)$ holds for some $i$, $j$, $k$ and
$R_{\ell} = 0$ for all $\ell \neq i$, $j$, $k$. 
\end{enumerate}
\end{theorem}
\noindent
This will be proved in section \ref{proofofT1}.  Next, 
we find when every other interior point of the cone has the unique
interior point of the appropriate Minkowski product of 
$P_{\tree}(\br)$ as a summand. 
Theorem \ref{T1} allows us to significantly narrow our search
for $P_{\tree}(\br)$ which satisfy this condition. 
In order to carry this out we bring in an alternative description
of the weightings $\omega \in P_{\tree}(\br)$, which can
be found in \cite{HMM} and \cite{HMSV}. 
We start by considering a weighting $\omega$
on a single trinode, $\tau$. Since $\Delta_2(\omega(\tau, 1), \omega(\tau, 2), \omega(\tau, 3))$
always holds, we may apply the 1-1 transformation of cones

$$T: P_3 \to \R_+^3$$
given by $T(\omega)(x_{ij}) = \frac{1}2(\omega(\tau, i) + \omega(\tau, j) - \omega(\tau, k)).$
This is an isomorphism of semigroups when $\R^3$ is given the standard lattice.  We represent
the image of $\omega$ via the $piping$ $model$ on the leaves of the tree, where the number of pipes
going from $i$ to $j$ is $T(\omega)(x_{ij})$, an example is picture below, see also figures \ref{Fig7} and \ref{Fig3}.

\begin{figure}[htbp]
\centering
\includegraphics[scale = 0.4]{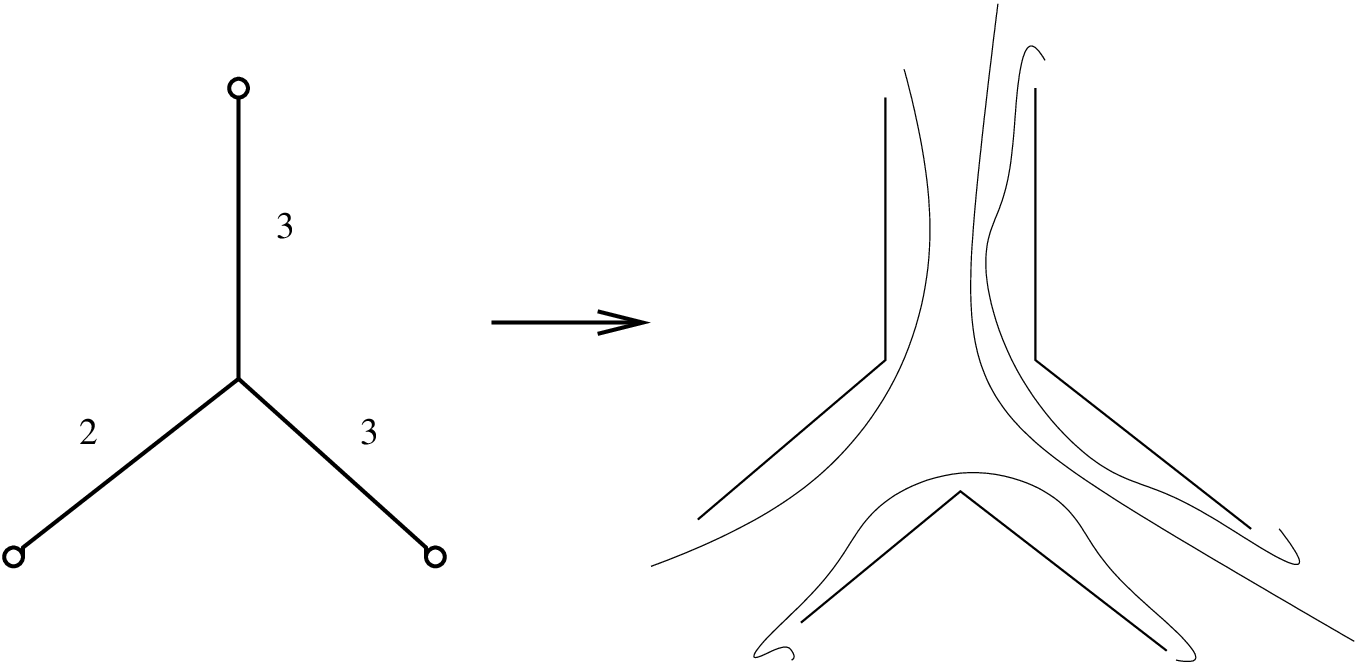}
\caption{Piping model for a trinode.}
\label{Fig6}
\end{figure}
\noindent
The map $T$ has 
an inverse $S$ given by $S(\eta)(\tau, i) = \eta(x_{ij}) + \eta(x_{ik})$. Using the map $T$ we may go from weightings on the trinode $\tau$ to planar graphs on the set 
$\{(\tau, i)\}$.
The transformation $T$ is useful as it clarifies
divisibility issues in the semigroup of lattice points in $P_3$: $\omega$ divides
$\omega'$ if and only if $T(\omega)(x_{ij}) \leq T(\omega')(x_{ij})$ for all $i,$ $j$. 
For a general tree $\tree$, something similar happens for graphs on the
set $L(\tree)$.  Given a graph $G$ on the set $L(\tree)$, we may construct a simultaneous
weighting of the edges of each trinode $\tau \in V(\tree)$ as follows.  
For each "pipe" $e \in G$ we consider the unique path $\gamma$ in $\tree$ joining the
end points of $e$, each trinode edge $(\tau, i)$ traversed by this path gets weight
$+1$.  We call this map $S_{\tree}$.  There is a section to this map, $T_{\tree}$, which is defined by
the following algorithm.  Given a weighting $\omega$ of $\tree$, consider
the simultaneous weighting of trinodes given by restricting $\omega$
to each $\tau \in V(\tree)$.  Apply $T$ to each of these weightings, and 
join up the ends of the resulting pipes in the unique way such that
the resulting graph is planar.  We leave it to the reader to show
that this is well-defined (that is, there is a unique planar graph for such a weighting,
this can also be found in \cite{HMM}.). 
An example is illustrated in figure \ref{Fig7}.\\

\begin{figure}[htbp]
\centering
\includegraphics[scale = 0.2]{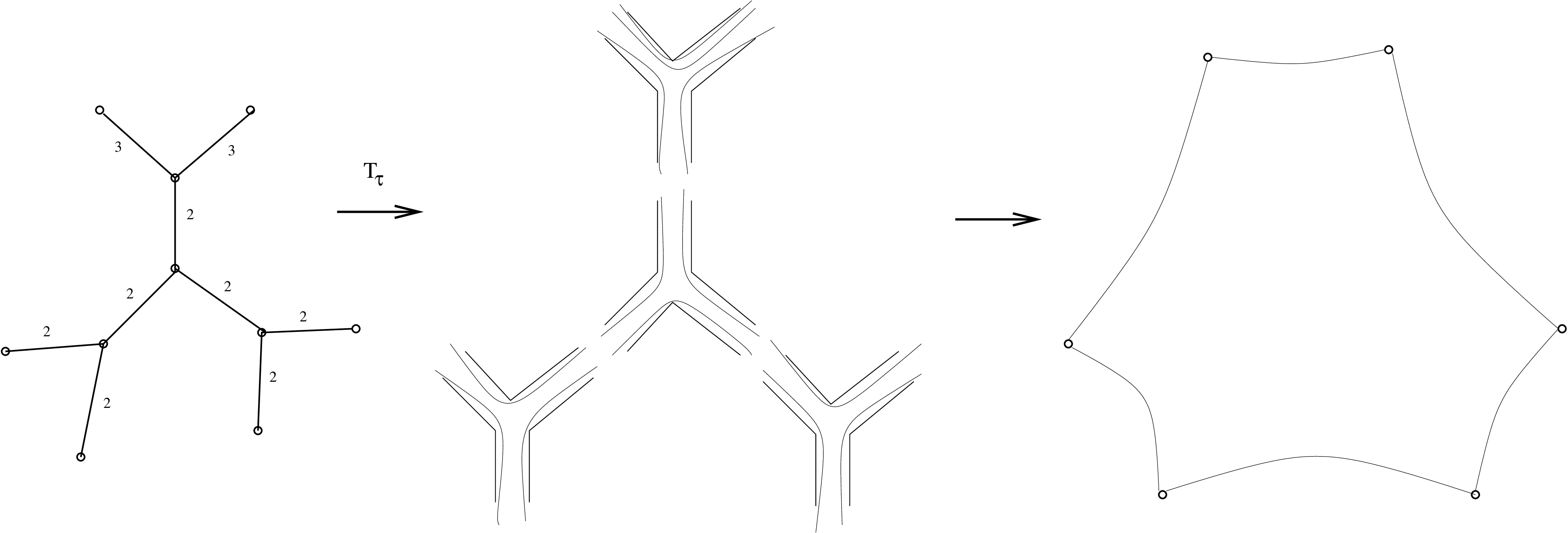}
\caption{Piping model for a general $\tree$.}
\label{Fig7}
\end{figure}

\noindent
The map $T_{\tree}$ cannot be an inverse as there are many (non-planar) graphs 
that give the same weighting under $S_{\tree}$, this redundancy accounts
for the elements in a certain basis of $\C[M_{\br}]$, linearly related
by the Pl\"ucker equations, which degenerate to the same element under the flat degeneration
defined by $\tree$. For more on this see \cite{HMM}. The reader will also note that a weighting
$\omega$ divides another weighting $\omega'$ if and only if 
$\omega|_{\tau}$ divides $\omega'|_{\tau}$ for all $\tau \in V(\tree)$, 
so the map $T$ is still very useful for questions of divisibility in the case of general $\tree$. 
Note that the piping model only makes sense for a tree $\tree'$ which
has a specified embedding into $\R^2$, because of our need to deal
with planar graphs.  For this reason we restrict our attention
to trees $\tree$ with such an embedding assumed.  This has no effect on our
results because the isomorphism class of $P_{\tree}$ as a multigraded
algebra is determined by the topological type of $\tree$, and therefore any $\C[S_{\tree}(\br)]$
is isomorphic to some $\C[S_{\tree'}(\br')]$ with $\tree'$ planar
and $\br',$ a permutation of the entries of $\br.$

We let $N_{ij}(\omega)$ be the multiplicity of edges between $i$ and $j$ in the graph $T_{\tree}(\omega)$.
We are now ready to state our main theorem. From now on
we denote the unique minimal degree internal weighting
for the cone on $P_{\tree}(\br)$ by $\omega_{\br}(\tree)$, if
it exists.  Also, we let $2_{\tree}$ be the weighting which
assigns $2$ to each edge of $\tree$, an example is illustrated
in figure \ref{Fig7} above.  In general, $T_{\tree}(2_{\tree})$
is always a complete planar cycle on the set $L(\tree)$. 

\begin{theorem}\label{maintheorem}
$\C[S_{\tree}(\br)]$ is Gorenstein if and only if
\begin{enumerate}
\item $a\br$ is as in theorem \ref{T1} for some $a$.
\item In this degree, $N_{ij}(\omega_{\br}(\tree) - 2_{\tree}) \geq n-4$ when it is nonzero.
\end{enumerate}
\end{theorem}
\noindent
This will be proved in section \ref{proofofmaintheorem}.   
We finish with a theorem which restricts the $a$-invariant
of $\C[S_{\tree}(\br)]$.

\begin{theorem}\label{a-inv}
$$a(\C[S_{\tree}(\br)]) \ | \ 2(n-2)$$
\end{theorem}
\noindent
This is proved in section \ref{ainv}. We mention here that this
theorem also restricts the $a$-invariant for $\C[M_{\br}]$ 
as well, because the $a$-invariant may be read off the Hilbert 
function of the algebra.

\section{Proof of theorem \ref{T1}.}\label{proofofT1}

In this section we will classify the polytopes
$P_{\tree}(\br)$ that have a unique interior point. 
This is the first step to understanding when the ring
$\C[S_{\tree}(\br)]$ is Gorenstein. Let $P_{\tree}$ 
denote the cone of weightings
$\omega$ on the tree $\tree$ such that

$$\Delta_2(\omega(\tau, 1), \omega(\tau, 2), \omega(\tau, 3))$$
holds for each internal vertex $\tau \in V(\tree)$.  A weighting
$\omega$ is on a face of this cone if and only if one of these
triangle inequalities is an equality for some vertex $\tau$.
Let $D_n$ be the cone of side lengths for $n$-sided
polygons (so $D_3 = P_3$), there is a map of cones $\pi: P_{\tree} \to D_n$ given by 
forgetting the weights on internal edges of $\tree$.  
We have the following identification, see prop 4.4 of \cite{HMSV}.

$$P_{\tree}(\br) = \pi^{-1}(\br).$$
We study the dimension of these fibers when $\br_i \neq 0$.
For $n=3$ there is nothing to say, so suppose $n > 3$. 
If for every side length $\br_i < \Sigma_{j \neq i} \br_j$ 
then we may form a planar $n$-gon
$\mathcal{P}$ with side lengths $\br$ with non-zero area, 
this in turn implies we may find such an $n$-gon 
where all triangles in $\mathcal{P}$ formed by the diagonals and sides
have non-zero area. Let $\omega$ be the (not necessarily integer) 
weighting of $\tree$ formed by the diagonal
lengths and side lengths of $\mathcal{P}$. 

\begin{figure}[htbp]
\centering
\includegraphics[scale = 0.35]{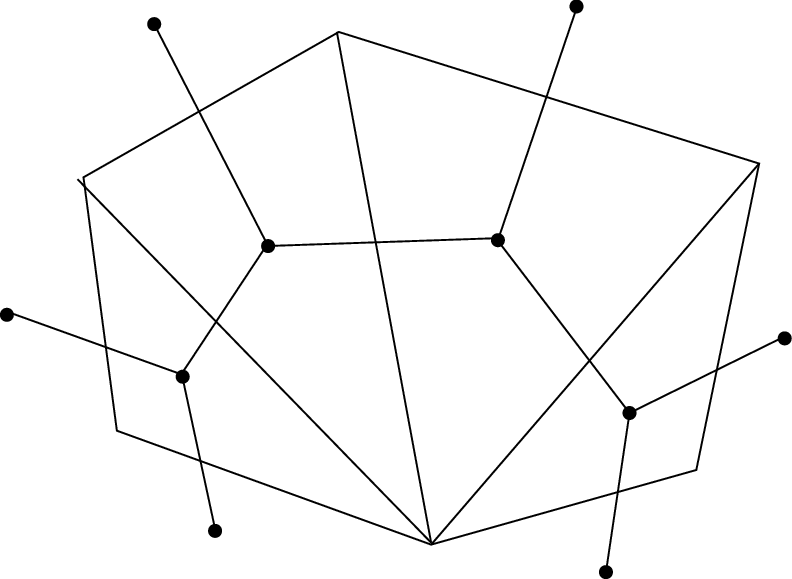}
\caption{The edges of $\tree$ recieve weight equal to the length of the dual diagonal in $\mathcal{P}$}
\label{tri}
\end{figure}

\begin{figure}[htbp]
\centering
\includegraphics[scale = 0.35]{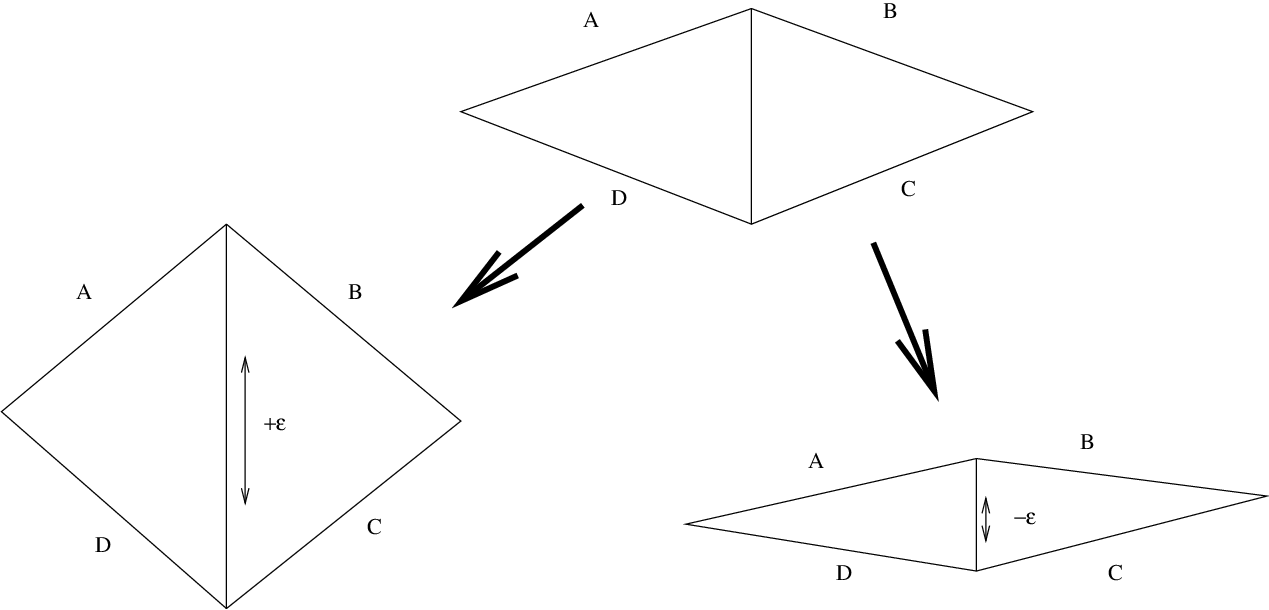}
\caption{Creating a neighborhood of $\omega$. The case
for a general polygon can be reduced to the case of
a quadrilateral by excising all but the $4$ edges
incident on the diagonal in question.}
\label{Fig1}
\end{figure}

Now, observe that for any diagonal $d$ of $\mathcal{P}$ specified by $\tree$ which borders
two triangles with non-zero area, we may stretch and contract the length of $d$ 
within some neighborhood $\epsilon$ without changing 
the lengths of any other sides and specified diagonals non-zero, 
see figure \ref{Fig1}.

\noindent
This implies that there is a small neighborhood of dimension $|I(\tree)|$ in $P_{\tree}(\br)$
which contains $\omega$. On the other hand, if some entry $\br_j$ 
of $\br$ has $\br_j = \sum_{i \neq j} \br_i$  then
$P_{\tree}(\br)$ can be nothing but a single point. Thus we conclude
that the dimension of $P_{\tree}(\br)$ can be either
$|I(\tree)|$ or $1.$ This discussion has some bearing
on the following proposition.

\begin{proposition}
Let $\br_i \neq  0$ for all $i.$
A weighting $\omega$  is in the interior of $P_{\tree}(\br) = \pi^{-1}(\br)$ only if it is in the interior of $P_{\tree}$. 
\end{proposition}

\begin{proof}
First note that if $\br$ has some entry equal to the sum of the other entries, 
then $P_{\tree}(\br)$ is a fiber over a point in a facet of $D_n$, 
therefore the unique point of $P_{\tree}(\br)$ is on a facet of $P_{\tree}$, so 
suppose this is not the case.  If $\omega$ is an interior point of $P_{\tree}(\br)$
then there is a $|I(\tree)|$-dimensional neighborhood of $\omega$
in $P_{\tree}(\br),$ this implies that the weighting $\omega(e)$ of each edge $e \in E(\tree)$
can be expanded and contracted while keeping the weight in $P_{\tree}(\br)$ without changing the other weights on the members of $L(\tree)$.  If some triangle inequality is an equality, then some
triangle $\tau$ in the polygon defined by $\omega$ is degenerate.  Let $e$ be the edge of $\tree$ dual to the longest edge in $\tau$, by the degeneracy of $\tau,$ the length of this edge cannot be increased without also increasing the lengths of one of the other two edges in $\tau.$  This implies that all triangle inequalities must be strict on $\omega,$ so this weight
is an internal point in $P_{\tree}.$
\end{proof}

\noindent
The content of this proposition can also be found in \cite{KM}. 
The following gives an algebraic characterization of the interior lattice points
of $P_{\tree}$ and $P_{\tree}(\br)$. 

\begin{proposition}\label{prop1}
A non-negative integer weighting $\omega \in P_{\tree}$ is in the interior
if and only if $\omega = \eta + 2_{\tree}$ for
some $\eta \in P_{\tree}$. 
\end{proposition}

\begin{proof}
If $\omega$ is in the interior of $P_{\tree}$ then
all inequalities defined by the condition $\Delta_2$
are strict.  After converting to the piping model,
we must have $T_{\tree}(\omega)(x_{ij}(\tau)) \geq 1$,
for each trinode $\tau \in V(\tree)$. 
This implies that $\omega$ has $2_{\tree}$ as a factor. 
Running this argument in reverse gives the converse. 
\end{proof}

\begin{corollary}\label{cor1}
If $\omega \in P_{\tree}(\br)$ is an interior
point, then $\omega = 2_{\tree} + \eta$ for
some $\eta \in P_{\tree}(\br - \vec{2})$
\end{corollary}

\begin{corollary}
The toric algebra $\C[P_{\tree}]$ is Gorenstein.
\end{corollary}
\noindent
This second corollary also follows from the same theory that 
gave us theorem \ref{deform}, and the fact that
the algebra of the Pl\"ucker embedding of $Gr_2(\C^n)$
is Gorenstein.  As a consequence of proposition \ref{prop1} and
corollary \ref{cor1} we get the following proposition. 

\begin{proposition}\label{prop2}
$P_{\tree}(\br)$ has a unique interior
point if and only if $\br = \vec{2} + \vec{R}$
such that $P_{\tree}(\vec{R})$ is a single point. 
\end{proposition}
\noindent
The next proposition classifies all 
$\vec{R}$ which have this property. 

\begin{proposition}\label{prop3}
$P_{\tree}(\vec{R})$ is exactly one point
if and only if $\vec{R}$ is of one of the following forms
\begin{enumerate}
\item $R_i = \Sigma_{i \neq j} R_j$ for some $i$. 
\item $\Delta_2(R_i, R_j, R_k)$ holds for some $i$, $j$, $k$ and
$R_{\ell} = 0$ for all $\ell \neq i$, $j$, $k$. 
\end{enumerate}
\end{proposition}

\begin{figure}[htbp]
\centering
\includegraphics[scale = 0.2]{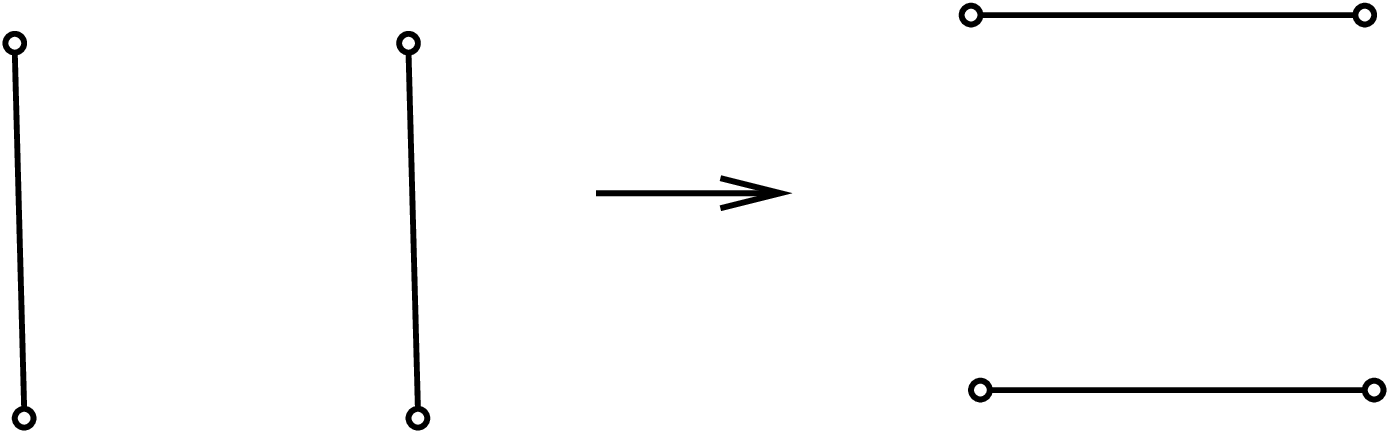}
\caption{Creating another weighting.}
\label{Fig3}
\end{figure}

\begin{proof}
For sufficiency, note that there is exactly
one polygon fitting the description given by
both cases above.  For necessity, we consider the
piping model $T_{\tree}(\omega)$ of the tree weighting. We suppose
$P_{\tree}(\br)$ is a single point and classify
the pipe arrangements allowed for weightings of $\tree$. 
Suppose we were allowed an arrangement of pipes
where two edges do not share any common vertex.  
Then we may swap pipes while maintaining the edge
weights, as in figure \ref{Fig3}. This implies that any pair of pipes must share a
common vertex.  There are exactly two ways for this
to happen, the reader can verify that 
$S_{\tree} \circ T_{\tree}(\omega) = \omega$
must satisfy the edge weight conditions in the
statement of the theorem. 
\end{proof}

\begin{corollary}\label{cor2}
$P_{\tree}(\br)$ has a unique interior
lattice point if and only if $P_{\tree'}(\br)$
has a unique interior lattice point, for all $\tree'$. 
\end{corollary}

When we convert the above weighting conditions to their graphical
representation on the set $L(\tree)$, we get the possibilities represented
below in Figure \ref{Fig4}.  One possibility is a graph where every pipe shares
a common incident vertex, the second possibility
has exactly three vertices with incident pipes.  
Propositions \ref{prop2} and \ref{prop3} then prove
theorem \ref{T1}.

\begin{figure}[htbp]
\centering
\includegraphics[scale = 0.33]{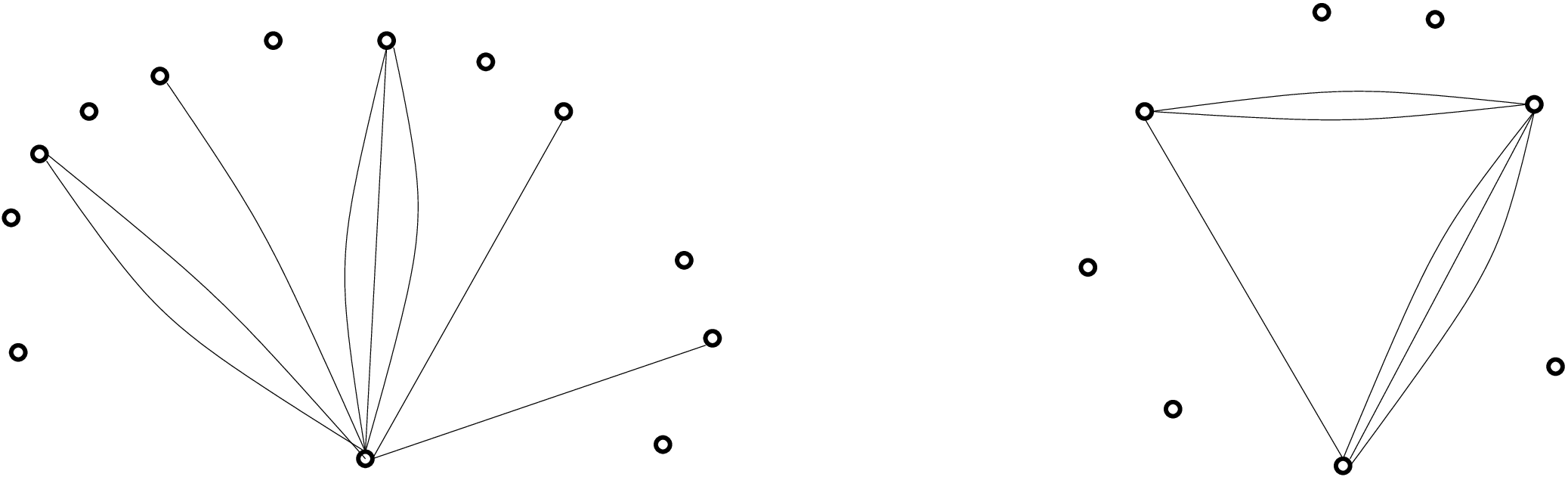}
\caption{Associated graphs for $\vec{R}$ from the proof of proposition \ref{prop3}.}
\label{Fig4}
\end{figure}

\section{Proof of theorem \ref{maintheorem}.}\label{proofofmaintheorem}
Theorem \ref{T1} gives a necessary condition
for $\C[S_{\tree}(\br)]$ to be Gorenstein.   
Now we see what must be added in order to ensure that all
interior lattice points carry the unique interior lattice
point $\omega_{\br}(\tree)$ as a summand.  
We will make use of the piping model
for most of this section.  For the cases presented in the statement
of theorem \ref{T1}, the first case has 
$N_{ij}(\omega_{\br}(\tree) - 2_{\tree}) = R_j$ and 
$N_{kj}(\omega_{\br}(\tree) - 2_{\tree}) = 0$ for
all $k$, $j \neq i$, and the second case has 
$N_{ij}(\omega_{\br}(\tree) - 2_{\tree}) = \frac{1}2(R_i + R_j - R_k)$ 
and $N_{m \ell}(\omega_{\br}(\tree) - 2_{\tree}) = 0$
for $\ell$ or $m \neq i$, $j$, $k$. 

\begin{proposition}
Let $\br = \vec{2} + \vec{R},$ where $\vec{R}$ satisfies the conditions of
proposition \ref{prop3}. 
$\C[S_{\tree}(\br)]$ is Gorenstein if and only
if there is no interior weighting $\omega$ in degree $k \geq a$
such that $N_{ij}(\omega - 2_{\tree}) < N_{ij}(\omega_{\br}(\tree) - 2_{\tree})$
for all $i, j.$
\end{proposition}

\begin{proof}
After converting $\omega$ to the piping model
and removing the complete cycle on $L(\tree)$ corresponding
to $2_{\tree}$ we get the graph of $\omega - 2_{\tree}.$  It is clear 
that if $N_{ij}(\omega - 2_{\tree}) \geq N_{ij}(\omega_{\br}(\tree) - 2_{\tree})$ for all $i,$ $j$ then $\omega_{\br}(\tree)$ is a summand of $\omega$. 

For the converse, suppose $N_{ij}(\omega -2_{\tree}) <  N_{ij}(\omega_{\br}(\tree) - 2_{\tree})$
for some $i, j \in L(\tree).$  We find a weighting $\omega'$ on a new tree $\tree'$ which has
a pair of leaves $i'$ and $j'$ connected to a common trinode $\tau$, with the number of
pipes between $i'$ and $j'$ in the trinode equal to $N_{ij}(\omega).$ To do this, simply exchange members of $L(\tree)$ with a permutation
$\sigma$ so that $\sigma(i)= i'$ and $\sigma(j) = j'$
are next to each other, and choose a $\tree'$ such that these leaves
are now incident on a common internal vertex $\tau.$ 
Carrying the graph corresponding to $\omega$ along with the permutation $\sigma$
produces a new graph which may have crossings, but this does not matter, as no crossings
can be introduced between $i'$ and $j'.$
We consider the weighting 
$S_{\tree'}\circ \sigma \circ T_{\tree}(\omega) = \omega'$. 
By corollary \ref{cor2} there exists a unique internal lattice point $\omega_{\sigma(\br)}(\tree')$
in the polytope $P_{\tree'}(\sigma(\br))$ with $N_{i'j'}(\omega_{\sigma(\br)}(\tree')) = N_{ij}(\omega_{\br}(\tree))$.  
By construction we have $N_{i'j'}(\omega' - 2_{\tree'}) < N_{i'j'}(\omega_{\sigma(\br)}(\tree')- 2_{\tree'})$, implying that
$\omega_{\sigma(\br)}(\tree')|_{\tau}$ cannot divide $\omega'|_{\tau}.$
This implies that $\omega_{\sigma(\br)}(\tree')$ cannot divide $\omega',$
and that $\C[S_{\tree'}(\sigma(\br))]$ is not Gorenstein.
The permutation group $\mathcal{S}_n$ acts on the algebra
of global sections of $Gr_2(\C^n)$ given by the Pl\"ucker embedding  
by permuting the entries of the multigrading, so we 
get $\C[M_{\br}] \cong \C[M_{\sigma(\br)}].$ Now 
by theorem \ref{deform} and corollary \ref{equiv},
$\C[S_{\tree}(\br)]$ cannot be Gorenstein either. 
\end{proof}

\begin{figure}[htbp]
\centering
\includegraphics[scale = 0.33]{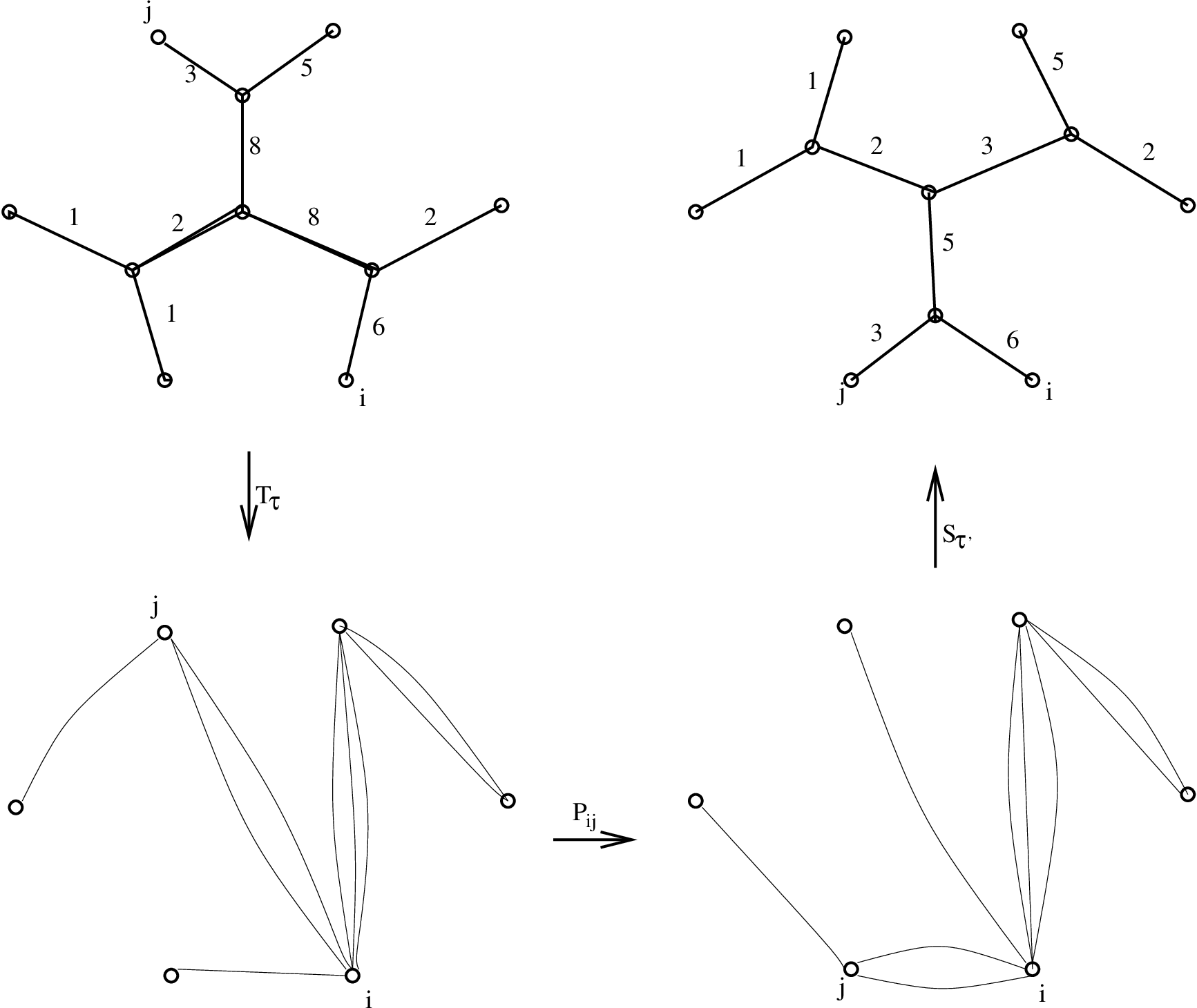}
\caption{Proof of theorem \ref{maintheorem}}
\label{Fig5}
\end{figure}

Now we are ready to prove theorem \ref{maintheorem}, 
this is accomplished with the next proposition. 

\begin{proposition}
For any $\C[S_{\tree}(\br)]$ such that some multiple
of $\br$ satisfies the criteria of theorem \ref{T1}, 
there is a weighting $\omega$ which has 
$N_{ij}(\omega - 2_{\tree}) < N_{ij}(\omega_{\br}(\tree) - 2_{\tree})$
if and only if $N_{ij}(\omega_{\br}(\tree) - 2_{\tree})$ is less than $n-4$
when it is nonzero. 
\end{proposition}

\begin{proof}
We must show that a weighting $\omega$ can be created with 
 $N_{ij}(\omega - 2_{\tree}) < N_{ij}(\omega_{\br}(\tree) - 2_{\tree})$
if and only if $N_{ij}(\omega_{\br}(\tree) - 2_{\tree})$ is less than $n-4.$
The inequality 
$N_{ij}(\omega -2_{\tree}) < N_{ij}(\omega_{\br}(\tree) - 2_{\tree})$
holds for some $\omega$ in degree $k$ only if we may start
with $N_{ij}(\omega - 2_{\tree})$  edges between $i$ and $j$
and add in edges to obtain a graph which is faithful to the condition that the edge
weights must be a multiple $k\vec{r}$.  It is necessary to have

$$\Sigma_{\ell \neq i, j} [\frac{k}a (R_{\ell} + 2) - 2] 
- [\frac{k}a(R_i + 2) - 2] - [\frac{k}a(R_j + 2) - 2] + 
2N_{ij}(\omega_{\br}(\tree) - 2_{\tree}) > 0,$$
where $a$ is the degree of $\omega_{\br}(\tree)$ and $k$ is the degree of $\omega$.  
To see this, note that the $\omega - 2_{\tree}$ weight on the the $\ell-$th leaf of $\tree$
must be

$$\frac{k}{a}a r_{\ell} - 2 = [\frac{k}{a}(R_{\ell} +2) - 2].$$

\noindent
It then follows that the quantity

$$[\frac{k}a(R_i + 2) - 2] + [\frac{k}a(R_j + 2) - 2]$$

\noindent
must be less than or equal to $2N_{ij}(\omega - 2_{\tree}) + \sum_{\ell \neq i, j} [\frac{k}{a}(R_{\ell} +2) - 2]$
Since we assumed $N_{ij}(\omega_{\vec{r}}(\tree) - 2_{\tree}) > N_{ij}(\omega - 2_{\tree}),$ we obtain
the above inequality.  

It remains to be seen how this inequality reduces to $N_{ij}(\omega_{\br}(\tree) - 2_{\tree}) < n-4$. 
In the case where $R_1 = \sum_{j \neq 1} R_j,$ we have $R_j = N_{1j}(\omega_{\tree}(\br) - 2_{\tree})),$ 
so the inequality reduces to 

$$2[\frac{k}{a} -1](n-4) > 2[\frac{k}{a} -1](N_{1j}(\omega_{\tree}(\br) - 2_{\tree}))$$

\noindent
for any of the non-zero nonzero $N_{1j}(\omega_{\tree}(\br) - 2_{\tree})$. This clearly
implies$N_{1j}(\omega_{\tree}(\br) - 2_{\tree}) < n- 4.$  Conversely, if $N_{1j}(\omega_{\tree}(\br) - 2_{\tree}) < n-4$
then we can recover this inequality for $k = 2a,$ and construct a graph $G$ with the desired properties as follows. 
We assume without loss of generality that $j = 2,$ and between the leaves $1$ and $2$ we put $R_2$ edges, note
that this is less than the required $R_2 + 1$ for $\omega_{\tree}(\br)$ to divide $\omega(G).$
Between $1$ and $j \neq 2$ we put $2R_j$ edges, and we add a complete planar cycle.  
The resulting graph requires $2$ more edges at each vertex $\neq 1, 2$ and $R_2 + 2$ more edges at vertices $1$ and $2$ to have the correct multi-degree. We have assumed that $R_2 < n-4,$ so it follows that $2R_2 + 4 < 2n - 4 = 2(n-2).$  This ensures
that there are enough spots left to assigne edges to $G$ in order to obtain the correct multidegree.

In the case where we have that $\Delta_2(R_1, R_2, R_3)$ holds with all other $R_{\ell} = 0,$ we may assume
without loss of generality that $a = 1,$ since all other $R_{\ell}$ are $0.$ Assuming $i = 2, j = 3,$
the above inequality becomes 

$$2k\frac{(R_1 - R_2 - R_3)}{2} + 2(k-1)(n-4) + 2N_{23}(\omega_{\tree}(\br) - 2_{\tree}) > 0 $$

\noindent
Using the identity $\frac{(R_1 - R_2 - R_3)}{2} = - N_{23}(\omega_{\tree}(\br) - 2_{\tree})$ this
becomes

$$2[k-1](n-4) > 2[k-1]N_{23}(\omega_{\tree}(\br) - 2_{\tree}).$$

\noindent
Since $k$ must be greater than $1,$ this implies the inequality.  Conversely, if
$N_{23}(\omega_{\tree}(\br) - 2_{\tree}) < n-4$ then we may recover the inequality
above for $k = 2,$ and construct a graph $H$ with the desired properties as follows.  
Between $2$ and $3$ we place $N_{23}(\omega_{\tree}(\br) - 2_{\tree})$ edges,
note that this is less than the number needed for $\omega(H)$ to carry $\omega_{\tree}(\br)$
as a divisor. We complete this to a cycle by adding a single edge between each consecutive pair 
$(3,4), \ldots (n-1, n), (n, 1), (1, 2).$  Now we add $2N_{12}(\omega_{\tree}(\br) - 2_{\tree}) + 1,$
and $2N_{12}(\omega_{\tree}(\br) - 2_{\tree}) + 1$ edges between $1, 2$ and $1, 3$ respectively. 
To finish, we must place edges in such a way that $2$ and $3$ each receive $N_{23}(\omega_{\tree}(\br) - 2_{\tree}) + 4 \leq n -1$
more edges.  There are $2(n-2)$ spots left to fill from the remaining vertices, so this is always possible.  
\end{proof}

\section{The $a$-invariant}\label{ainv}
Since the polytopes $P_{\tree}(\br)$ are the fibers
of $\pi$, a morphism of convex cones induced
by ambient linear map, we get $P_{\tree}(k\br) = kP_{\tree}(\br)$.
This allows us to prove theorem \ref{a-inv}.  This
theorem is implied by the following proposition.

\begin{proposition}
If $P_{\tree}(k\br)$ has a unique internal
lattice point then $k$ must divide $2(n-2)$
\end{proposition}

\begin{proof}
If $n \leq 3$ then our ring is isomorphic to $\C[x]$. 
Furthermore, if any $k\br_i = 2$ then $k = 1$ or $2$.
This takes care of all cases except when $R_i = \Sigma_{i \neq j} R_j$
and all $R_j > 0$.  Note that $k$ must divide

$$\Sigma_{i \neq j} (R_j + 2) - (R_i + 2) = R_i + 2(n-1) - R_i -  2 = 2(n-2)$$
\end{proof}

\begin{example}[Gorenstein property first shown by B. Howard and M. Herring, \cite{HH}]
Consider the case $\br = (1, \ldots, 1) = \vec{1}$,
this case satisfies all the conditions
of theorem $\ref{maintheorem},$ with the unique interior
point occurring in the polytope $P_{\tree}(\vec{2})$, 
the lattice points of which give the degree $2$ part of the algebra.  
Therefore $\C[S_{\tree}(\vec{1})]$ and $\C[M_{\vec{1}}]$ are
Gorenstein, with $a$-invariant equal to $2.$  The latter algebra
is of particular importance in \cite{HMSV}.
\end{example}

\begin{example}
In order to see the range of possible $a$-invariants, we'll look 
at a small example.  Consider the weights $(1, 1, 2, 4, 6),$  the third
graded component of $\C[M_{(1, 1, 2, 4, 6)}]$ has weights
$(3, 3, 6, 12, 18) = (2, 2, 2, 2, 2) + (1, 1, 4, 10, 16).$  Since
$16 = 1 + 1 + 4 + 10,$ and each number is greater than or equal 
to $5 - 4 = 1,$ this algebra is Gorenstein with the generator of the
canonical module in degree $3.$
\end{example}

\section{Acknowledgements}
I would like to thank my adviser John Millson
for many useful conversations about the Gorenstein property
in the case of the family $\C[M_{\br}]$. I'd also like to thank
Ben Howard for suggesting this problem, explaining to me the Gorenstein property
for toric algebras in terms of convex cones, and clarifying 
theorem \ref{deform}.

\end{document}